\documentstyle[12pt,amstex,amssymb]{amsart}
\newtheorem{th}{Theorem}

\newtheorem{cor}[th]{Corollary}
\newtheorem{lemma}[th]{Lemma}

\begin{document}
\title[Relative entropy for MASA]
{Relative entropy for maximal abelian subalgebras of matrices and 
the entropy of 
unistochastic matrices} 
\author[M. Choda]{Marie Choda}
\address{Osaka Kyoiku University, 
           Asahigaoka, Kashiwara 582-8582, Japan }

\thanks{Partially support by JSPS Grant No.17540193}

\date{}
\dedicatory{Dedicated to the memory of Professor Masahiro Nakamura}
\maketitle

\newcommand{\dom}{{\cal D}}
%
\begin{abstract} 
Let $A$ and $B$ be two  maximal abelian *-subalgebras of the $n\times n$ 
complex  matrices $M_n(\mathbb{C}).$  
To study the movement of the 
inner automorphisms of $M_n(\mathbb{C}),$
we modify the Connes-St$\o$rmer relative entropy 
$H(A | B)$ and the Connes relative entropy $H_\phi(A | B)$ 
with respect to a state $\phi,$  and introduce the two 
kinds of the constant 
$h(A | B)$ and $h_\phi(A | B).$ 
For the unistochastic matrix $b(u)$ 
defined by a unitary $u$ with $B = uAu^*,$  
we show that $h(A | B)$ is the entropy $H(b(u))$ of $b(u).$ 
This is obtained by our computation of $h_\phi(A | B).$ 
The $h(A | B)$ attains to the maximal value $\log n$ if and only if 
the pair $\{A, B\}$ is orthogonal in the sense of Popa. 
\end{abstract}
\section{Introduction} 
In a step to introduce the notion of the entropy for automorphisms on 
operator algebras,  
Connes and St$\o$rmer defined   in \cite{CS} 
the  relative entropy $H(A | B)$  
for finite dimensional von Neumann subalgebras 
$A$ and $B$ of a finite von Neumann algebra $M$ with a fixed normal normalized trace $\tau.$ 
\vskip 0.3cm

After then, the study about the  $H(A | B)$ is extended to two directions. 

One development was started  by Pimsner and Popa in \cite{PP} as 
the relative entropy $H(A | B)$   
for arbitrary von Neumann subalgebras $A$ and $B$ of $M,$ 
(see \cite{NS} for many interesting results in this direction). 

The other was a generalization due to Connes in \cite{C} 
by changing the trace $\tau$ to a  state  $\phi$  on $M.$ 
He defined  the relative entropy $H_\phi (A | B)$ 
for two subalgebras $A$ and $B$ of $M$ 
with respect to  $\phi.$  
\vskip 0.3cm

We modify the Connes-St$\o$rmer relative entropy 
$H(A | B)$ and the Connes relative entropy $H_\phi(A | B)$ 
with respect to a state $\phi,$  and we introduce the 
corresponding two kinds of the constant 
$h(A | B)$ and $h_\phi(A | B).$ 
In the case where $A = M,$ 
$h( M \mid B) $ is  nothing else but the Connes-St$\o$mer relative entropy 
$H( A \mid B).$  In general, $0 \leq h_\phi(A | B) \le H_\phi(A | B),$ 
and if $M$ is  an abelian von Neumann algebra, then 
$h_\phi( A \mid B) = H_\phi( A \mid B)$ so that it 
coincides with the conditional entropy in the ergodic theory. 

In this paper, we restrict our subjects to the maximal abelian subalgebras 
(abbreviated as MASA's)  of the $n\times n$ complex  matrices  
$M_n(\mathbb{C}) $. 

If $A$ and $B$ are two MASA's of $M_n(\mathbb{C}),$ then 
there exists a unitary matrix 
$u$ with $B = uAu^*,$ which we denote by  $u(A, B).$ 
Each unitary matrix $u$ induces a unistochastic matrix $b(u),$ 
which is a typical example of a bistochastic matrix. 

For a bistochastic matrix $b,$ it is introduced in \cite{ZSKS} 
the notion of the entropy $H(b)$ and the weighted entropy 
$H_\lambda(b)$ with respect to $\lambda.$ 

Here we show a relation between $h(A \mid B)$  
for a pair $\{A, B\}$ of MASA's  and the $H(b(u(A,B))).$ 
This is a consequence of our computation on 
$h_\phi(D \mid uDu^*),$  
where $\phi$ is a state of $M_n(\mathbb{C})$ and 
$D$ is the set of $n \times n$ diagonal matrices 
corresponding to the eigenvectors of $\phi.$  
The $h_\phi(D \mid uDu^*)$ satisfies the following relation :   
$$h_\phi(D \mid uDu^*) 
= H_\lambda(b(u)^*) + S(\phi|_D) - S(\phi|_{uDu^*}) ,$$
where $\lambda = \{\lambda_1, \cdots, \lambda_n\}$ is  
the eigenvalues for $\phi$ and $S(\psi)$ means the entropy of 
a positive linear functional $\psi.$ 

As the special case where $\phi$ is the normalized trace, we have 
$$h(A \mid B) =  H(b(u(A,B))).$$ 
\smallskip

In the case of the Connes-St$\o$rmer relative entropy, 
for a given two MASA's $A$ and $B$ of $M_n(\mathbb{C}),$ 
$H(A \mid B)$ is not equal to $H(b(u(A,B)))$ in general. 
For an example, see  \cite [Appendix] {PSW}.

\smallskip

The above results show that $h_\phi(A | B)$ for MASA's of $M_n(\mathbb{C})$ 
are determined by the entropy for the related unistochastic matrices 
in the sense of \cite{ZSKS}. Also, the value $h(D \mid uDu^*) $ 
is   depending on the inner automorphism 
$Ad_u$ defined by the unitary $u,$ and we can consider these values as
a kind of conditional entropy for $Ad_u$ conditioned by $D.$ 
\smallskip

Two  MASA's $A$ and $B$ of $M_n(\mathbb{C}) $ 
are {\it orthogonal} in the sense of Popa \cite{Po} if 
$$A \cap B = {\mathbb{C}}1 \quad \text{and} \quad 
E_A E_B = E_B E_A = E_{A\cap B} = E_{{\mathbb{C}}1},$$ 
where 
$E_A$ is the conditional expectation of $M_n(\mathbb{C})$ onto $A$ 
such that $\tau \circ E_A = \tau$ for the normalized trace $\tau$  of 
$M_n(\mathbb{C})$. 
This means   that 
\begin{eqnarray*}
\lefteqn{\ \quad \ B \ \subset \ M_n(\mathbb{C}) }\\
&& \cup  \ \ \ \quad \ \cup  \\
&& {\mathbb{C}} 1 \  \subset \  A 
\end{eqnarray*} 
is a commuting square in  the sense of \cite{GHJ}. 
\vskip 0.3cm

We have that a pair $\{A, B\}$ of MASA's is orthogonal 
if and only if $h(A | B)$ takes the maximal value, and then the value 
is $\log n.$ 

\bigbreak
\section{Preliminaries}
In this section, we summarize notations, terminologies and 
basic facts which we need later. 

Let $M$ be a finite von Neumann algebra, and let $\tau$ be a fixed 
faithful normal tracial state. 
By a von Neumann subalgebra $A$ of $M,$ we  
mean that $A$ has the same identity with $M.$ 
Let $E_A$ be the $\tau$-conditional expectations on $A.$ 

\subsection{Relative  entropy }  
\subsubsection{Relative   entropy of Connes-St$\o$rmer } 
First, we review about the formulae of the relative  entropy 
by Connes and St$\o$rmer in \cite{CS} (cf. \cite{NS}). 

Let $S$ be the set of all finite families $(x_i)$ of positive 
elements in $M$ with $1 = \sum_i x_i.$ 
Let $A$ and $B$ be two von Neumann subalgebras of $M.$ 

The  relative  entropy $H (A | B)$   is 
$$H( A \mid B) = 
\sup_{(x_i) \in S} \sum_{i} (\tau \eta E_{B}(x_i) - \tau \eta E_A(x_i)). $$
Here $\eta(t) = -t \log t, (0< t\leq 1)$ and $\eta(0) = 0.$ 

Let $\phi$ be a normal state on $M.$ 
Let $\Phi$ be the set of all finite families $(\phi_i)$ of positive 
linear functionals on $M$ with $\phi = \sum_i \phi_i.$ 

The relative  entropy $H_\phi(A | B)$ of $A$ and $B$ with respect to $\phi$ 
is 
$$ H_\phi (A | B) = 
\sup_{(\phi_i) \in \Phi} 
\sum_i ( S(\phi_i\mid_A, \phi \mid_A) - S(\phi_i\mid_B, \phi \mid_B) )$$ 
where $S( \psi| \varphi)$ is the relative entropy for 
two positive linear functionals $\psi$ and $\varphi,$ 
and  $H_\tau(A | B) = H( A \mid B).$ 

\subsubsection{Relative   entropy of positive linear functionals } 
After, we need the precise form $S( \psi| \varphi)$ 
in the case of finite dimensional algebras $C,$  
mainly the full matrix algebra $M_n(\mathbb{C})$ 
(the set of the $n\times n$  matrix $x = (x(i,j))_{ij}$ with $x(i,j) 
\in \mathbb{C}$). 

We denote by Tr the canonical trace on $C,$ 
that is, Tr($p) = 1$ for every minimal projection $p \in C.$ 

Let $\psi$ be a positive linear functional on $C.$ 
We denote by $Q_\psi$ the density operator of $\psi,$ 
that is,  $Q_\psi \in C$ is a unique positive operator with 
$$\psi(x) = \text{Tr}(Q_\psi x), \quad (x \in C),$$ 
and the von Neumann  entropy of $\psi$ is given by 
$$S(\psi) = \text{Tr}(\eta(Q_\psi)),$$ 

Let $\psi$ and $ \varphi $ be two positive linear functionals on $C.$ 
If $\psi \leq \lambda \varphi $ for some $\lambda > 0,$ then 
the relative  entropy of $\psi$ and $\varphi $ is given as 
$$S(\psi, \varphi) = \text{Tr}(Q_\psi (\log Q_\psi - \log Q_\varphi )),$$  
(cf. \cite{NS}, \cite{OP}). 
\smallskip

\vskip 0.3cm

\subsubsection{Conditional relative entropy} 
Let $A$ and $B$ be von Neumann subalgebras of $M.$ 
We modify $H(A | B)$ and  $H_\phi(A | B)$  and  
as a replacement of $H(A | B)$ (resp. $H_\phi(A | B)$ ) 
we define the following constant $h(A | B)$ (resp. $h_\phi(A | B)$). 
\smallskip 
 
The {\it conditional relative entropy} $ h( A \mid B)$ 
of  $A$ and $B$  conditioned by $A$ is 
$$ h( A \mid B) 
= \sup_{(x_i) \in S}  
\sum_{i} (\tau \eta E_{B}(E_{A}(x_i)) - \tau \eta E_{A} (x_i)). $$

Let $S(A) \subset S$ be the set of all finite families $(x_i)$ of positive 
elements in $A$ with $1 = \sum_i x_i.$ 
Then it is clear that 
$$ h( A \mid B) 
= \sup_{(x_i) \in S(A)}  
\sum_{i} (\tau \eta E_{B}(x_i) - \tau \eta (x_i)). $$

Let $S'(A) \subset S(A)$ be the set of all finite families $(x_i)$ 
with each $x_i$ a scalar multiple of a projection in $A.$ 
Then by the same proof with in  \cite{PP}, 
$$ h( A \mid B) 
= \sup_{(x_i) \in S'(A)}  
\sum_{i} (\tau \eta E_{B}(x_i) - \tau \eta (x_i)). $$ 
Hence, if $A$ is finite dimensional, 
we only need to consider the families consisting of  
scalar multiples of orthogonal minimal projections. 
\vskip 0.3cm 
Let $\phi$ be a normal state of $M.$ 

The {\it conditional relative entropy} of $A$ and $B$ 
with respect to $\phi$ conditioned by $A$ is 
$$ h_\phi (A | B) = 
\sup_{(\phi_i) \in \Phi} 
 \sum_i ( S(\phi_i\mid_A, \phi \mid_A) - 
S((\phi_i\circ E_A)\mid_B, (\phi\circ E_A)\mid_B) ).$$ 

Assume that $M$ is finite dimensional and that 
the density matrix of $\phi$ is contained in $A.$ 
Let $\Phi(A) \subset \Phi$ be the set of all finite families 
$(\phi_i)$ whose density operators $(Q_i)_i$ are contained  in  $A.$ 
Since the density matrix of $\phi_i \circ E_A$ is $E_A(Q_i),$ we have 
$$ h_\phi (A | B) = 
\sup_{(\phi_i) \in \Phi(A)} 
\sum_i ( S(\phi_i\mid_A, \phi \mid_A) 
   - S(\phi_i \mid_B, \phi \mid_B) ).$$ 
\vskip 0.3cm

\subsubsection{\bf Remark.} 
We will study about $ h_\phi (A | B)$ and $h( A \mid B)$ 
for von Neumann subalgebras $A$ and $B$ 
of a general finite von Neumann algebra $M$ elsewhere. 
Here we just remark the following facts. 

(1) If $\phi$ is the normalized trace $\tau,$ then by 
\cite[ Theorem 2.3.1(x) ] {NS} 
$$ h_\tau (A | B) = h( A \mid B).$$ 

(2) It is clear that 
$h( M \mid B) $ is  nothing else but the Connes-St$\o$mer relative entropy 
$H( M \mid B),$ and we have the relation with Index 
by Pimsner-Popa \cite{PP}. 

(3) In general, $0 \leq h_\phi(A | B) \le H_\phi(A | B).$ 
If $M$ is  an abelian von Neumann algebra, then 
$h_\phi( A \mid B) $ coincides with the conditional entropy 
in the ergodic theory 
by a similar proof as in \cite [ p. 158 ] {NS}.
\vskip 0.3cm

\noindent
\subsection{Unistochastic matrices and the entropy} 
When a matrix $x \in M_n(\mathbb{C})$ is given, 
we denote the $(i,j)$-component of $x$ by $x(i,j).$ 
A  matrix  $b \in M_n(\mathbb{C})$ is called {\it bistochastic} if 
$b(i,j) \geq 0$ for all $i,j = 1, \cdots, n,$ 
$\sum_{i=1}^n b(i,j) = 1 $ for all $j = 1, \cdots, n$ 
and $\sum_{j=1}^n b(i,j) = 1 $ for all $i = 1, \cdots, n.$ 
Let $\lambda = \{\lambda_1, \cdots, \lambda_n \}$ be a probability vector.  
\smallskip

The weighted entropy $ H_\lambda(b)$ of a bistochastic matrix $b$ 
with respect to $\lambda $ and 
the entropy $H(b)$ for a bistochastic matrix $b$  are given  in \cite{ZSKS} 
as the following forms, respectively : 
$$  H_\lambda(b) = \sum_{k=1}^n \lambda_k \sum_{j=1}^n \eta(b(j,k)), $$
and 
$$ H(b) =  \frac 1n \sum_{i=1}^n \sum_{j=1}^n \eta ( b(i,j) ).$$ 
\smallskip

Let $u$ be a $n\times n$ unitary matrix. 
The {\it unistochastic matrix $b$ defined by $u$} is the bistochastic matrix 
given as 
$$b(i,j)= |u(i,j)|^2, \ (i, j = 1,2, \cdots, n).$$   

\bigbreak
\section{Results} 

\begin{lemma} 
Let $A$ be a maximal abelian subalgebra of $M_n(\mathbb{C}),$ and let 
$\{p_1, \cdots, p_n\}$ be the minimal projections of $A.$  

(1) If $\psi$ is a positive linear functional of $M_n(\mathbb{C}),$ 
then 
$$ S(\psi|_A) = \sum_{j=1}^n \eta( \psi(p_j) ).$$ 

(2) If $\phi$ is a state of $M_n(\mathbb{C})$ and if 
$\phi = \sum_i \phi_i$ is a finite decomposition of $\phi$ 
into a sum of positive linear functionals, 
then 
$$ \sum_i S(\phi_i|_A, \phi|_A) = - \sum_i S(\phi_i|_A) + S(\phi|_A).$$ 
\end{lemma}

\begin{pf} 

(1) Since the density operator of $\psi|_A$ is written as 
$\sum_j \psi(p_j) p_j,$ we have 
$$ S(\psi|_A) 
= \text{Tr}_A(\sum_j \eta( \psi(p_j) ) p_j ) 
= \sum_{j=1}^n \eta( \psi(p_j) ).$$

(2) We denote the density operator of $\phi_i|_A$ by $Q_i$ and 
the density operator of $\phi|_A$ by $Q$, then 

\begin{eqnarray*} 
\lefteqn{\sum_i S(\phi_i \mid_A, \phi \mid_A) = \sum_i \text{Tr} (Q_i(\log Q_i - \log Q) ) } \\
&&= - \sum_i S(\phi_i \mid_A) - \text{Tr}( \sum_i Q_i \log Q) \\
&&= - \sum_i S(\phi_i \mid_A) + S(\phi \mid_A).
\end{eqnarray*}
\end{pf} 

Let $\phi$ be a state of $M_n(\mathbb{C}).$ 
We  number the set of the eigenvalues of  $Q_\phi$  as 
$\lambda = \{\lambda_1 \geq \lambda_2 \geq \cdots \geq \lambda_n\}.$ 
Let us  decompose $Q_\phi $ into the form   
$Q_\phi = \sum_{i=1}^n \lambda_i e_i$, where $\{e_1, \cdots, e_n\}$ 
is a family of mutually orthogonal minimal projections 
in $M_n(\mathbb{C}),$ which we fix. 
Let $D$ be the MASA generated 
by the projections $\{e_1, \cdots, e_n\},$ which we denote by $D(\phi)$ 
when we need. 

Let $\{e_{kl} \}_{k,l = 1, \cdots, n}$ be a system of matrix units of 
$M_n(\mathbb{C})$ 
such that $e_{ii} = e_i$ for all $i = 1, \cdots, n.$ 
We give the matrix representation for each $x \in M_n(\mathbb{C})$ 
depending on this 
matrix units $\{e_{kl} \}_{k,l = 1, \cdots, n}$ 
so that $D = D(\phi)$ is the diagonal  algebra. 
Let $u  \in M_n(\mathbb{C})$ be a  unitary,  
and let $b(u)$ be the unistochstic matrix defined by $u.$ 
Under these situations, we have the following theorem.  
\vskip 0.3cm

\begin{th} 
Let $\phi$ be a state of $M_n(\mathbb{C}),$ and 
let $D$ be the diagonal algebra of $M_n(\mathbb{C}).$ 
Let $u \in M_n(\mathbb{C})$ be a unitary, then 
$$h_\phi(D \mid uDu^*) 
= H_\lambda(b(u)^*) + S(\phi|_D) - S(\phi|_{uDu^*}). $$
\end{th}
\vskip 0.3cm

\begin{pf}
First we remark that $Q_\phi$ is contaied in $D$ and 
$S(\phi) = S(\phi|_D).$ 
Using the matrix representation of the unitary $u,$ 
the matrix representation of each $ue_ju^*$ is given by 
$$ue_ju^* = \sum_{k,l} u(k,j) \overline{u(l,j)}e_{k,l} 
= (u(k,j)\overline{u(l,j)})_{kl}.$$
Let $(\phi_i)_{i = 1, \cdots, n}$ be the positive linear functionals  of 
$M_n(\mathbb{C})$ such that $Q_{\phi_i} = \lambda_i e_i$  for all $i.$ 
Then $\phi = \sum_i \phi_i$ gives a  finite decomposition of $\phi,$ 
and 
\begin{eqnarray*}
\lefteqn{
\sum_j\eta(\phi_i(ue_ju^*)) = \sum_j\eta(\lambda_i \mid u(i,j) \mid^2)  }\\
&&= \sum_j ( \ \eta(\lambda_i )\mid u(i,j) \mid^2 
    + \lambda_i \eta(\mid u(i,j) \mid^2) \ ) \\
&&= \eta(\lambda_i ) + \lambda_i \sum_j \eta(\mid u(i,j) \mid^2). 
\end{eqnarray*} 
Hence by Lemma 1, 
\begin{eqnarray*}
\lefteqn{
\sum_{i =1}^n ( S(\phi_i\mid_D, \phi \mid_D) - S(\phi_i\mid_{uDu^*}, \phi \mid_{uDu^*} ) ) }\\
&&= - \sum_i \sum_j\eta(\phi_i(e_j)) + \sum_j\eta(\lambda_j)  \\
&& \quad  +  \sum_i \sum_j\eta(\phi_i(ue_ju^*))   -  \sum_j \eta(\phi( ue_ju^*)) \\
&&= \sum_i \sum_j\eta(\phi_i(ue_ju^*))   -  \sum_j \eta(\phi( ue_ju^*)) \\
&&= S(\phi|_D)  + H_\lambda(b(u)^*)  - S(\phi|_{uDu^*}).
\end{eqnarray*} 
Thus 
$$H_\phi(D \mid uDu^*) 
\geq   H_\lambda(b(u)^*) + S(\phi|_D) - S(\phi|_{uDu^*}).$$
\vskip 0.3cm

To prove the opposite inequality,  assume that 
$(\phi)_{i \in I}$ be a given finite family which is contained 
in $\Phi(D).$ 
We denote $Q_{\phi_i}$ simply by $Q_i.$ 

Then  
$$\lambda_l = \sum_{i \in I} Q_i(l,l), \quad \text{and} \quad 
  Q_i(l,k) = 0 \quad \text{if} \quad l \ne k,$$ 
for all $l, k = 1, \cdots, n$ and $i \in I,$ 
and 
$$\phi_i(ue_ju^*) =   \sum_k Q_i(k, k) |u(k,j)|^2 
\quad \text{and} \quad 
\phi_i(e_j) =  Q_i(j, j),$$ 
for all $i \in I$ and $j = 1, \cdots, n.$ 

Since  
$\eta(st) = \eta(s)t + s\eta(t)$ and 
$\eta(s + t) \le \eta(s) + \eta(t)$ for all positive numbers $s$ and $t,$ 
by using that  $\sum_j|u(k,j)|^2 = 1$ for all 
$k,$  we have that 
\begin{eqnarray*}
\lefteqn{ 
\sum_i S(\phi_i|_{uDu^*})  -\sum_i S(\phi_i|_D)} \\
&&= \sum_i \sum_j \eta(  \sum_k Q_i(k,k) | u(k,j) |^2) 
  - \sum_i \sum_j \eta(Q_i(j,j))  \\
&&\leq  
\sum_i \sum_j \{\sum_k \eta(Q_i(k,k)) |u(k,j) |^2 
+ \sum_k Q_i(k,k) \eta(|u(k,j) |^2) \}\\
&& \quad \quad 
  - \sum_i \sum_j \eta(Q_i(j,j)) \\
&& = \sum_i \sum_k \eta(Q_i(k,k)) (\sum_j |u(k,j) |^2 ) 
 + \sum_j \sum_k \{(\sum_i Q_i(k,k)) \eta(|u(k,j) |^2)\}  \\
&& \quad \quad    - \sum_i \sum_j \eta(Q_i(j,j)) \\
&& =  \sum_i \sum_k \eta(Q_i(k,k)) 
 + \sum_j \sum_k \lambda_k \eta(|u(k,j) |^2) 
 - \sum_i \sum_j \eta(Q_i(j,j)) \\
&& = H_\lambda(b(u)^*)
\end{eqnarray*}

These imply that  by Lemma 1,   
\begin{eqnarray*}
\lefteqn{
\sum_{i \in I} \{ S(\phi_i\mid_D, \phi \mid_D) - 
     S(\phi_i\mid_{uDu^*}, \phi \mid_{uDu^*} )  \} }\\
&&= -\sum_i S(\phi_i|_D) + S(\phi|_D) 
    + \sum_i S(\phi_i|_{uDu^*}) - S(\phi|_{uDu^*}) \\
&& \leq   H_\lambda(b(u)^*) + S(\phi|_D)  - S(\phi|_{uDu^*}) 
\end{eqnarray*} 

Hence we have always that 
$$h_\phi(D|uDu^*) 
\leq  H_\lambda(b(u)^*) + S(\phi|_D) - S(\phi |_{uDu^*}),$$
and 
$$h_\phi(D|uDu^*) 
= H_\lambda(b(u)^*) + S(\phi|_D) - S(\phi |_{uDu^*}).$$
\end{pf}

Let $D$ be the algebra of diagonal $n\times n$ matrices. 
Let $\phi$ be a state of $M_n(\mathbb{C})$, and let 
$v$ be a unitary such that $vDv^* = D(\phi).$ 
Let $\phi_v$ be the state given by 
the inner perturbation by $v$ : $\phi_v(x) = vxv^*.$ 
Then we have the following relation between $h(A | B)$ and the set of 
$\{h_\phi(A | B)\}_\phi.$

\vskip 0.3cm
\begin{cor} 
Let $D$ be the algebra of diagonal $n\times n$ matrices, and let 
$u \in M_n(\mathbb{C})$ be a unitary. 
Then 
$$
h(D \mid uDu^*)  = H(b(u)) =  \frac 1n \sum_{i=1}^n \sum_{j=1}^n 
  \eta ( \mid u(i,j) \mid^2 ) $$ 
$$ \quad = \max_\phi h_{\phi_v} (D \mid uDu^*), $$
where the maximum is taken over all states $\phi$ of $M_n(\mathbb{C}).$
\end{cor}

\begin{pf} 
The first equality is clear because 
$h(D \mid uDu^*) =  h_\tau(D \mid uDu^*) $ 
and 
$H_\lambda(b(u)^*) = 
\frac 1n \sum_{i=1}^n \sum_{j=1}^n \eta ( \mid u(i,j) \mid^2 )$ 
where the eigenvalues for $\tau$ is $\lambda = \{1/n, \cdots, 1/n\}.$ 
We denote the minimal projections of $D$ by $\{e_1, \cdots, e_n\}.$ 
Let $\phi$ be a state, and  let 
$\{\lambda_1, \lambda_2, \cdots, \lambda_n\}$ be eigenvalues of 
$Q_\phi.$ 
We decompose $Q_\phi = \sum_{i=1}^n \lambda_i p_i$ 
by mutually orthogonal minimal projections $\{p_1, p_2, \cdots, p_n \}.$ 
Let $v$ be a unitary with $ve_iv^* = p_i$ for all $i.$
Then    
$$Q_{\phi_v} = \sum_i \lambda_i e_i = {v^*} Q_\phi v.$$ 
That is, 
$D(\phi_v) = D$ and $S(\phi_v) = S((\phi_v)|_D).$  
Since 
$$S(\phi_v) \leq S((\phi_v)|_{uDu^*})$$ 
(cf. \cite[Theorem 2.2.2 (vii) ] {NS}), 
we have by Theorem 2 that,  
for the unistochastic matrix $b(u)$ defined by $u,$ 
\begin{eqnarray*}
\lefteqn{h_{\phi_v} (D \mid uDu^*)}\\
&&=H_\lambda(b(u)^*) + S((\phi_v)|_D) - S((\phi_v)|_{uDu^*}) \\
&&\leq H_\lambda(b(u)^*)  \\
&&\leq H(b(u)^*) = H(b(u)) = h(D \mid uDu^*).
\end{eqnarray*}
\end{pf}

\vskip 0.3cm
As another easy consequence of Theorem 2,  by 
the property of the normalized trace $\tau,$ we have 
the following statement for general MASA's $A$ and $B$  :  

\begin{cor} 
Let $A$ and $B$ be maximal abelian subalgebras of $M_n(\mathbb{C}).$ 
Then 
$$h(A \mid B) =  H(b(u)) = \frac 1n \sum_{i=1}^n \sum_{j=1}^n \eta ( \mid u(i,j) \mid^2 )$$
where  the $(u(i,j))_{ij}$ for $u = u(A,B)$  
is given by the matrix units whose minimal  projections generates $A.$ 
\end{cor}

\vskip 0.3cm
\begin{cor} 
Let $\{A_0, B_0\}$ be a pair of maximal abelian subalgebras of $M_n(\mathbb{C}).$ 
Then $\{A_0, B_0\}$ is an orthogonal pair if and only if 
$$h(A_0 \mid B_0) = \log n = \max \ h(A \mid B),$$
where the maximum is taken over 
the set of a pair $\{A,B\}$ of maximal abelian subalgebras of $M_n(\mathbb{C}).$ 
\end{cor}

\begin{pf}
Let  $u \in M_n(\mathbb{C})$ be a unitary with $uAu^* = B.$  
By a characterization of Popa(\cite{Po}), 
a pair $\{A, B\}$ of MASA's is orthogonal if and only if 
$\tau(aubu^*) = 0 $ for all $a,b \in A$ with $\tau(a) = \tau(b) = 0.$ 
This implies that $\{A, B\}$ is orthogonal if and only if 
$$|u(j,k)| =  1 / {\sqrt n} \quad \text{for all} \quad j, k.$$ 
Hence we have the conclusion.
\end{pf}


\end{document}